\documentclass[11pt]{article}
\usepackage{amsmath,amssymb,amsthm}

\newcommand{\Tbar}{\ensuremath{\overline{\mathcal{T}}}}
\newcommand{\Mbar}{\overline{\mathcal{M}}}
\newcommand{\hess}{\operatorname{Hess}}
\newcommand{\grad}{\operatorname{grad}}
\newcommand{\gradl}{\operatorname{grad}\ell}
\newcommand{\basis}{\{\lambda_{\alpha},J\lambda_{\alpha},\gradl_{\beta}\}_{\alpha\in\sigma,\beta\in\tau}}
\newcommand{\tr}{\operatorname{tr}}
\newcommand{\dt}{\frac{d}{dt}}
\newcommand{\dr}{\frac{\partial}{\partial r}}
\newcommand{\dthet}{\frac{\partial}{\partial \vartheta}}
\newcommand{\mathR}{\mathbb{R}}

\newcommand{\caB}{\mathcal{B}}
\newcommand{\caF}{\mathcal{F}}

\newcommand{\caP}{\mathcal{P}}

\newcommand{\caT}{\mathcal{T}}

\newtheorem{definition}{Definition}
\newtheorem{lemma}[definition]{Lemma}
\newtheorem{theorem}[definition]{Theorem}
\newtheorem{corollary}[definition]{Corollary}
\begin{document}

\title{Understanding Weil-Petersson curvature\footnote{2000 Mathematics Subject Classification Primary: 32G15; Secondary: 20H10, 30F60.}}     
\author{Scott A. Wolpert}
\date{September 7, 2008}         
\maketitle

\begin{center}
{\em Dedicated to Shing-Tung Yau}
\end{center}

\begin{abstract}
A brief history of the investigation of the Weil-Petersson curvature and a summary of Teichm\"{u}ller theory are provided.  A report is presented on the program to describe an intrinsic geometry with the Weil-Petersson metric and geodesic-length functions.  Formulas for the metric, covariant derivative and formulas for the curvature tensor are presented.  A discussion of methods is included.  Recent and new applications are sketched, including results from the work of Liu-Sun-Yau, an examination of the Yamada model metric and a description of Jacobi fields along geodesics to the boundary.
\end{abstract}

\section{Introduction}       

This written lecture is in honor of the birthday of Shing-Tung Yau, a most distinguished mathematician and teacher.  The present year is the anniversary of Weil's notes on using the Uniformization theorem to introduce an $L^2$ metric for the space of Riemann surfaces, \cite[see articles 1958b, c]{Wlwk2}.  

Recent progress on finite-dimensional Weil-Petersson (WP) geometry is discussed in the introductions of \cite{LSY1, LSY2, Wlcomp, Wlbhv} and in the notes \cite{LSY5, Wlpers, Wlgeom}.  A brief history of the investigation of WP curvature is as follows.  Ahlfors established formulas for the connection and curvature in terms of singular integrals and found that holomorphic sectional curvatures are negative \cite{Ahcurv}.  Royden showed that the holomorphic sectional curvature is bounded above by $-1/2\pi(g-1)$ \cite{Royicm}, followed then by the results of Tromba \cite{Trcurv} and the author \cite{Wlchern} including the general negativity of all sectional curvatures.  The curvature tensor formula of \cite{Wlchern} was generalized by Siu \cite{Siuwp} and Schumacher \cite{Schu} to the setting of deformation spaces of K\"{a}hler-Einstein metrics.  Schumacher considered bounds and growth for curvature, and established strong negativity of the curvature in the sense of Siu \cite{Schucurv}.  Huang and Teo continue the examination of the sectional curvatures of the metric \cite{Zh,Zh4,Zh2,Zh3,LPT}.  Liu, Sun and Yau as part of their program of investigation of the canonical metrics on Teichm\"{u}ller space show that the metric is dual Nakano negative \cite{LSY3}.

In \cite{Wlchern, Wlhyp} the approach was also applied to investigate the family hyperbolic metric for the vertical line bundle of the universal curve over the moduli space of Riemann surfaces.  The Chern form computed from the hyperbolic metric and considered as a current was found to represent the expected Chern class in cohomology, even though the form is singular at infinity \cite{Wlhyp}. The Chern form for the renormalized family hyperbolic metric for the tangent lines to cusps was calculated in \cite{Wlcusps}.   Freixas generalizes and uses these results in his investigations and applications of arithmetic intersection theory \cite{GFM}.  Trapani showed that the WP Ricci form represents in cohomology for the Deligne-Mumford compactification the determinant of the logarithmic tangent bundle, even though the Ricci form is singular at infinity \cite{Trap}.  Liu, Sun and Yau show that all the Chern forms for the WP metric considered as currents represent the expected Chern classes on the Deligne-Mumford compactification \cite{LSY3,LSY4}.  Applications are presented. 

We now report on the program to describe an intrinsic geometry with the WP metric and geodesic-length functions.  An invariant of a hyperbolic metric is the length $\ell_{\alpha}$ of the closed geodesic $\alpha$ in a free homotopy class.  Elements of the intrinsic geometry are given by the following expansions: \\
for the Fenchel-Nielsen twist and almost complex structure
\[
2t_{\alpha}=J\gradl_{\alpha},
\]
for the root-length gradients $\lambda_{\alpha}=\grad \ell_{\alpha}^{1/2}$ and disjoint geodesics
\[
\langle\lambda_{\alpha},\lambda_{\alpha'}\rangle\ =\ (2\pi)^{-1}\delta_{\alpha\alpha'}\ +\ O((\ell_{\alpha}\ell_{\alpha'})^{3/2}) ,
\]
for the Fenchel-Nielsen twist parameter and a maximal collection of disjoint geodesics (a pants decomposition)
\[
\frac12\sum_{\alpha}\,d\ell_{\alpha}\wedge d\tau_{\alpha},
\]
for the covariant derivative
\[
D_U\lambda_{\alpha}=3\ell_{\alpha}^{-1/2}\langle J\lambda_{\alpha},U\rangle J\lambda_{\alpha}\ +\ O(\ell_{\alpha}^{3/2}),
\]
for the curvature tensor
\[
R(\lambda_{\alpha},\lambda_{\alpha},\lambda_{\alpha},\lambda_{\alpha})= \frac{3}{16\pi^3\ell_{\alpha}}\ +\ O(\ell_{\alpha}),
\]
and for disjoint geodesics not all the same
\[
R(\lambda_{\alpha},\lambda_{\beta},\lambda_{\gamma},\lambda_{\delta}) = O((\ell_{\alpha}\ell_{\beta}\ell_{\gamma}\ell_{\delta})^{1/2}),
\]
where all remainder term constants are uniform for bounded geodesic-lengths.  Beyond the single leading term contributions to the curvature tensor are bounded.  
We present and discuss the formulas in Section 3 and consider applications in Section 5. 

The above expansions reveal the {\em asymptotic factorization} of the WP metric.  In a region of Teichm\"{u}ller space where the geodesics $\sigma=\{\alpha\}$ have small length the metric is almost the product of the metrics for the complex lines $\{\lambda_{\alpha},J\lambda_{\alpha}\}$ and the tangential metric for the boundary space $\caT(\sigma)\,=\,\{\ell_{\alpha}=0,\alpha\in\sigma\}$ in the augmented Teichm\"{u}ller space $\Tbar$.  The complex lines $\{\lambda_{\alpha},J\lambda_{\alpha}\}$ define an extension of the normal bundle of the boundary space.  The boundary space is a product of lower dimensional Teichm\"{u}ller spaces with WP metrics.  Asymptotic factorization partially appears in Masur's original expansion \cite{Msext}.  We find the factorization behavior is at least $C^2$; factorization appears in the metric, connection and curvature tensor.  The metric for the complex lines also compares to Yamada's model metric $\pi^3(4dr^2\,+\,r^6d\theta)$.  The WP metric and Yamada model $C^2$-compare with a 
factor $1+O(r^4)$.  

The lecture is organized as follows.  Basic material on WP geometry and the augmented Teichm\"{u}ller space is presented in Section 2.  The intrinsic geometry of geodesic-lengths and the WP metric is presented in Section 3.  An argument for bounding holomorphic quadratic differentials on collars and a method for bounding the kernel of the operator $(D\,-\,2)^{-1}$ are described in Section 4.  In Section 5.1 select results of Liu, Sun and Yau are described.  In Section 5.2 we examine the comparison to the Yamada model and describe Clairaut equations for a WP geodesic and a small geodesic-length.  In Section 5.3 we consider the distance to a proper stratum of $\Tbar$ and examine projecting geodesics, the length minimizing paths to a stratum of $\Tbar$.  We discuss variation of a projecting geodesic for top dimensional strata and characterize the  Jacobi fields describing variation through projecting geodesics.  The first variation of distance to a stratum is given by the cosine formula of \cite{BH}.  We present the formula for the second variation of distance.

\section{Basics of Teichm\"{u}ller theory}
Let $\caT$ be the Teichm\"{u}ller space for homotopy marked genus $g$, $n$ punctured Riemann surfaces $R$ of negative Euler characteristic.  A point of $\caT$ is the equivalence class of a pair $(R,f)$, with $f$ a homeomorphism from a reference topological surface $F$ to $R$.  By Uniformization a conformal structure determines a unique complete compatible hyperbolic metric $ds^2$ for $R$.  The Teichm\"{u}ller space is a complex manifold of dimension $3g-3+n$ with the cotangent space at $R$ represented by $Q(R)$, the space of holomorphic quadratic differentials on $R$ with at most simple poles at the punctures.  Weil introduced a Hermitian cometric. The Riemannian pairing is the real part of the dual form.
\begin{definition}\index{Weil-Petersson metric!definition}
The Weil-Petersson cometric is $\langle\varphi,\psi\rangle\ = \ \int_R\varphi\overline\psi\,(ds^2)^{-1}.$
\end{definition}
In the 1940's Teichm\"{u}ller introduced the Finsler metric with conorm given as $\|\varphi\|_T=\int_R|\varphi|$.  \index{Teichm\"{u}ller metric!definition}
The Weil-Petersson (WP) dual metric is invariant under the action of the {\em mapping class group}, $\mathrm{MCG}$, \index{mapping class group} the group of orientation preserving homeomorphisms of $F$ modulo the subgroup of homeomorphisms homotopic to the identity relative to punctures. The WP metric projects to the quotient $\mathcal M=\caT/\mathrm{MCG}$, the moduli space of Riemann surfaces.  First properties are that the metric is K\"{a}hler, non complete with negative sectional curvature $\kappa$ with $\sup_\caT\kappa=0$ (except for $\dim\caT=1$ where 
$\sup_\caT\kappa<0$) and $\inf_\caT\kappa=-\infty$.  The metric continues to be the first metric for understanding the K\"{a}hler geometry of Teichm\"{u}ller space \cite{LSY1, LSY2, McM, Wlpos, Yun1}. {\em In practice and experience the WP geometry of $\caT$ corresponds to the hyperbolic geometry of surfaces.} 

A hyperbolic surface has a {\em thick-thin} \index{thick-thin decomposition} decomposition with {\em thin} the region of injectivity radius below a threshold value.  The thin components of a hyperbolic surface are neighborhoods of cusps or are collars (fixed area tubular neighborhoods of short geodesics).  Mumford first observed that the set of hyperbolic surfaces with lengths of closed geodesics bounded below by a constant $c>0$ forms a compact subset $\mathcal M_c$ of the moduli space $\mathcal M$. In general the totality of all thick regions of hyperbolic surfaces of a given topological type forms a compact set of metric spaces in the Gromov-Hausdorff topology.  The Teichm\"{u}ller and WP metrics are comparable on $\mathcal M_c$.  

We have for the reference topological surface $F$.
\begin{definition} \index{complex of curves}
A $k$-simplex of the complex of curves $C(F)$ is a set of $k+1$ distinct free homotopy classes of non trivial, non peripheral, mutually disjoint simple closed curves of $F$. The pants graph $\mathrm{P(F)}$ has vertices the maximal simplices 
of $C(F)$. \index{pants graph} Vertices of $\mathrm{P(F)}$ are connected by an edge provided the corresponding sets of free homotopy classes differ by replacing a single curve by a curve intersecting the original curve one or two times.  
\end{definition} 
The vertices of $C(F)$ are the free homotopy classes of non trivial, non peripheral simple closed curves. The simplices of $C(F)$ are the convex sums of vertices.  The complex of curves $C(F)$ is a lattice, partially ordered by inclusion of simplices with maximal simplices, called {\em pants decompositions}, having dimension $3g-4+n$.   A pants decomposition decomposes a surface into a union of 
$2g-2+n$ three holed spheres.  The pants graph $\mathrm{P(F)}$ becomes a metric space by specifying the edges to have unit-length. Corresponding to a pants decomposition are global coordinates, Fenchel-Nielsen (FN) coordinates for $\caT$ given as gluing-parameters for constructing surfaces from right hyperbolic hexagons \cite{Abbook}. The construction begins with right hexagons which can be doubled across alternating edges to obtain a pair of pants, a genus zero hyperbolic surface with three geodesic boundaries with lengths free-parameters in $\mathR_{> 0}$.  Boundaries of pants of a common length can be abutted to construct a larger surface.  

A pants decomposition for $F$ provides a combinatorial scheme to abut boundaries of pants to obtain a hyperbolic surface of genus $g$ with $n$ punctures.  In abutting boundaries there is the free-parameter of the relative displacement of one boundary with respect to the other.  Overall for each abutting there are two free-parameters.  The first parameter is the common boundary {\em geodesic-length} $\ell$ valued in $\mathR_{>0}$. The second parameter is the relative displacement $\tau$ valued in $\mathR$ measured in hyperbolic distance ($\tau$ is initially measured between appropriate footpoints and then analytically continued).
\begin{theorem}
For a pants decomposition $\caP$ the FN coordinates $\Pi_{j\in\caP}(\ell_j,\tau_j):\rightarrow(\mathR_{>0}\times\mathR)^{3g-3+n}$ provide a real analytic equivalence for $\caT$.  
\end{theorem}
 
The {\em augmented Teichm\"{u}ller space} $\Tbar$ is a partial compactification in the style of Bailey-Borel \cite{Abbook,HrCh}.  The space $\Tbar$ is important to the understanding of the Deligne-Mumford compactification of the moduli space $\mathcal M$ and for the WP geometry.  Frontier spaces are adjoined to $\caT$ corresponding to allowing geodesic-lengths $\ell_j$ to assume the value zero with the FN angles $\theta_j=2\pi\tau_j/\ell_j$ undefined (in polar coordinates the angle is undefined at the origin).  The vanishing $\ell_j=0$ describes a degenerate hexagon with a side-length vanishing and the adjacent sides meeting at a common point on the circle at infinity for the hyperbolic plane.  The vanishing $\ell_j=0$ corresponds to a degenerate hyperbolic surface with a pair of cusps in place of a simple closed curve ($\gamma_j$ in the pants decomposition is now represented by the horocycles about the cusps).  In general for a simplex $\sigma\subset C(F)$, the $\sigma$-null stratum is the space of structures $\caT(\sigma)=\{R\ \mbox{degenerate}\mid\ell_{\alpha}(R)=0\mbox{ iff }\alpha\in\sigma\}$. \index{complex of curves} The frontier spaces $\caF_{\caP}=\cup_{\sigma\subset\caP}\caT(\sigma)$ subordinate to a pants decomposition $\caP$ are adjoined to $\caT$ with a neighborhood basis for $\caT\cup\caF_{\caP}$ prescribed by the specification that
\[
((\ell_{\beta},\theta_{\beta}),\ell_{\alpha}):\caT\cup\caT(\sigma)\rightarrow\Pi_{\beta\notin\sigma}(\mathR_{>0}\times\mathR)\times\Pi_{\alpha\in\sigma}(\mathR_{\ge0})
\]
is a continuous map. For a simplex $\sigma\subset\caP,\,\caP'$, contained in two pants decompositions, the neighborhood systems are equivalent.  

A structural property is that the deformation spaces $\caT(\sigma)$ are products of lower dimensional Teichm\"{u}ller spaces and the limit of the tangential component of the WP metric of $\caT$ is simply the WP metric of the lower dimensional stratum $\caT(\sigma)$.  The role of the augmented Teichm\"{u}ller space in WP geometry first appeared in Masur's work \cite{Msext}. A sequence of marked hyperbolic surfaces $\{R_n\}$ converges in $\Tbar$ provided there is a simplex $\sigma$ contained in a pants decomposition $\mathcal P$ with $\ell_{\alpha}(R_n),\,\alpha\in\sigma$ limiting to zero (no convergence condition is placed on $\theta_{\alpha},\,\alpha\in\sigma$) and FN parameters $(\ell_{\beta},\theta_{\beta}),\,\beta\in\caP-\sigma$ converging. 
\begin{theorem} \index{augmented Teichm\"{u}ller space}
The augmented Teichm\"{u}ller space $\Tbar=\caT\cup_{\sigma\in C(F)}\caT(\sigma)$ is a non locally compact stratified space.  The augmented Teichm\"{u}ller space is the WP completion of Teichm\"{u}ller space.
\end{theorem}
A point of $\Tbar-\caT$ represents a marked degenerate hyperbolic structure for which a simplex of $C(F)$ has each element represented by a pair of cusps. \index{complex of curves} The augmentation is also described as the Chabauty topology closure of the discrete faithful cofinite representations of $\pi_1(R)$ into $\mathrm{PSL}(2;\mathR)$ modulo conjugacies by $\mathrm{PSL}(2;\mathR)$ \cite{HrCh}.  A sequence of discrete faithful representations converges to a union of conjugacy classes of discrete faithful representations of smaller topological type.  The quotient $\Tbar/\mathrm{MCG}$ is topologically the Deligne-Mumford stable curve compactification $\mathcal M$ of the moduli space of curves. 

The basic properties of the compactification of the moduli space and of the augmented Teichm\"{u}ller space are described in the original research announcement of Bers \cite{Bersdeg}. Definitions include a Riemann surface with nodes, holomorphic $q$-differentials (push downs of sections of powers of the relative dualizing sheaf of a stable curve) and a description of the Deligne-Mumford compactification, as well as descriptions of the augmented Teichm\"{u}ller space.  Bers included a description of fundamental sets for the action of the mapping class group.  For a pants decomposition $\caP$ let $MCG(\caP)$ be the subgroup of mapping classes stabilizing the maximal simplex $\caP$.  The subgroup of Dehn twists about elements of $\caP$ is finite index in $MCG(\caP)$.  

\begin{theorem}  For a constant $c_{g,n}$ depending on the topological type and for a pants decomposition $\caP$, the Bers region $\caB(\caP)=\{\ell_j< c_{g,n},j\in\caP\}\subset\caT$ is a fundamental set for $MCG/MCG(\caP)$ in particular $\caT\subset\cup_{h\in MCG/MCG(\caP)}h(\caB(\caP))$.  The region $\overline{\caB(\caP)}=\{0\le \ell_j< c_{g,n}\}\subset\Tbar$ is a fundamental set for $MCG/MCG(\caP)$ acting on $\Tbar$.  
\end{theorem} 

The pants graph $\mathrm{P(F)}$ with edge-metric provides a quasi isometric model for the WP metric on $\caT$ and $\Tbar$ as follows.  The Bers regions cover Teichm\"{u}ller space and satisfy $h(\caB(\caP))=\caB(h(\caP))$ for $h\in MCG$.  The bounded valence multivalued relation $ R\in\caT\leftrightarrow\{\caP\mid R\in\caB(\caP)\}$ is the basis for (non unique) embeddings $h:\caT \rightarrow \mathrm{P(F)}$ and $k:\mathrm{P(F)}\rightarrow\caT$.  A point is associated with the pants decomposition of a Bers region containing the point.  Brock established that the embeddings are quasi isometries (satisfy $d(x,y)/c'-c''\le d'(f(x),f(y))\le c'd(x,y)+c''$ for positive constants).  Brock and in joint work with Margalit established the following \cite{BrMr, Brkwp}. 
\begin{theorem}
$\caT$ and $\Tbar$ are quasi isometric to $\mathrm{P(F)}$.  For the special topological types $(1,2)$ and $(0,4)$ the spaces $\caT$ and $\Tbar$ are quasi isometric to the Farey graph.
\end{theorem}   

The augmented Teichm\"{u}ller space with the WP metric is a $\mathrm{CAT}(0)$ metric space \cite{DW2, DW3, Wlcomp, Yam2} (a complete, simply connected, generalized non positively curved 
space \cite{BH}). In particular $(\Tbar,d_{WP})$ is a length space, a metric space with unique distance-realizing paths (geodesics) between pairs of points.   Furthermore the WP metric has the Euclidean triangle comparison property: for a triangle in $\Tbar$ and a corresponding triangle in $\mathbb E^2$ with corresponding edge lengths, chords of the former (located by endpoints on sides) have lengths bounded by lengths of corresponding chords of the latter.  (The conditions $\mathrm{CAT}(0)$ and Gromov hyperbolicity are independent, even though the latter is a general negative curvature condition.)  Each closed stratum of the augmented Teichm\"{u}ller space is the closure of a single $\sigma$-null stratum.  The $\sigma$-null strata $\caT(\sigma)\subset\Tbar,\,\sigma\in C(F)$ are intrinsic to the WP metric geometry and $\Tbar$ is an infinite polyhedron as follows \cite{Wlcomp}. 

\begin{theorem} \index{augmented Teichm\"{u}ller space}
$\Tbar$ is a stratified metric space with each $\sigma$-null stratum, $\sigma\in C(F)$, characterized as the union of all open geodesics containing a given point of the stratum as an interior point.  The interior of a geodesic in $\Tbar$ is contained in a single $\sigma$-null stratum (geodesics do not refract at strata).  $\Tbar$ itself is characterized as the closed convex hull of the maximally degenerate hyperbolic structures (the unions of thrice punctured spheres).
\end{theorem}

The structure of strata provides that the extended $\mathrm{MCG}$ (both orientation preserving and reversing classes) is the full group of isometries of $\caT$ as follows \cite{MW,Wlcomp}.    A WP isometry extends to the completion $\Tbar$ and preserves the intrinsic strata structure, as well as the partial ordering of inclusion of simplices of $C(F)$.  Ivanov established that order preserving bijections of $C(F)$ are induced by elements of the extended $\mathrm{MCG}$ \cite{Ivaut}. \index{complex of curves} In particular for an orientation preserving isometry there is a corresponding element of $\mathrm{MCG}$ and the two maps agree on the maximally degenerate structures in $\Tbar$ (a maximally degenerate structure is uniquely determined by its simplex).  The two maps agree on $\Tbar$, the closed convex hull of the maximally degenerate structures. The isometry coincides with the element of $\mathrm{MCG}$.

\section{WP intrinsic geometry}

The correspondence between WP geometry and the hyperbolic geometry of surfaces begins with geodesic-length functions.  Associated to each non trivial, non peripheral free homotopy class $\alpha$ on a marked hyperbolic surface is the length $\ell_{\alpha}$ of the unique geodesic in the free homotopy class.   In the Uniformization geodesic-lengths are explicit with $2 \cosh \ell_{\alpha}/2=\tr A$ for the free homotopy class $\alpha$ corresponding to the conjugacy class of $A$ in the deck group within $\mathrm{PSL}(2;\mathR)$. 

Geodesic-length functions  have a direct relationship to WP geometry.  For a simple closed geodesic $\alpha$ on a hyperbolic surface, the surface can be {\em cut} along the geodesic to form two circle boundaries, which can then be identified by a relative rotation to form a new surface.  A Hamiltonian flow on $\caT$ is defined by considering families of surfaces $\{R_t\}$ for which at time $t$ reference points from sides of the original geodesic are displaced by $t$ on $R_t$.  Collections of geodesic-length functions provide local coordinates for Teichm\"{u}ller space.   We describe the WP metric, connection and curvature tensor in terms of gradients of geodesic-length functions.
\begin{definition}
For the free homotopy class $\alpha$ of a closed curve let: $\ell_{\alpha}$ be the geodesic-length function, $\gradl_{\alpha}$ the associated WP gradient, $\lambda_{\alpha}=\gradl^{1/2}_{\alpha}$ the root-length gradient, and $t_{\alpha}$ the infinitesimal FN twist vector field.
\end{definition}

A basic feature of WP geometry is the twist-length duality $2t_{\alpha}=J\gradl_{\alpha}$ for $J$ the almost complex structure of $\caT$ \cite{WlFN}.  For a simplex $\sigma\subset C(F)$ consider geodesic-length functions in a neighborhood of an augmentation point $p$ of $\caT(\sigma)\subset\Tbar$.  We refer to the elements of $\sigma$ as the {\em short geodesics.}

\begin{definition}
A relative length basis for a point $p$ of $\caT(\sigma)$ is a collection $\tau$ of vertices of $C(F)$ disjoint from the elements of $\sigma$ such that at $p$ the gradients $\{\gradl_{\beta}\}_{\beta\in\tau}$ provide the germ of a frame over $\mathbb R$ for the tangent space $T\caT(\sigma)$.
\end{definition}  

The basic expansion for the WP metric and length gradients is the following \cite{Wlbhv}.

\begin{theorem}
\label{wppair}
The pairing of geodesic-length gradients for disjoint geodesics $\alpha,\beta$ satisfies
\[
0\ <\ \langle\gradl_{\alpha},\gradl_{\beta}\rangle\ -\ \frac{2}{\pi}\ell_{\alpha}\delta_{\alpha\beta}\quad \mbox{is}\quad O(\ell^2_{\alpha}\ell^2_{\beta})
\]
where for $c_0$ positive the remainder term constant is uniform for $\ell_{\alpha},\,\ell_{\beta}\le\,c_0$.
\end{theorem}

The expansion is based on Riera's exact formula for gradient pairings and on estimating double coset sums of $e^{-2dist}$ for the uniformization group $\Gamma\subset PSL(2;\mathbb R)$.  For a point $p$ of $\caT(\sigma)$ and relative length basis $\tau$ the tangent fields $\{\lambda_{\alpha}, J\lambda_{\alpha}, \gradl_{\beta}\}$ provide a local frame for $T\caT$.  A refined counterpart \cite{Wlbhv} of Masur's original expansion is a consequence.  For short geodesics $\alpha\in\sigma$ and elements $\beta\in\tau$ of a relative length basis
\begin{gather*}
\langle\lambda_{\alpha},\lambda_{\alpha'}\rangle\ =\ (2\pi)^{-1}\delta_{\alpha\alpha'}\ +\ O(\ell_{\alpha}^{3/2}\ell_{\alpha'}^{3/2}) \\
\langle \lambda_{\alpha},\gradl_{\beta}\rangle\ =\ O(\ell^{3/2}_{\alpha})
\end{gather*}
where by twist-length duality and disjointness
\[
\langle J\lambda_{\alpha},\lambda_{\alpha'}\rangle\ =\ \langle J\lambda_{\alpha},\gradl_{\beta}\rangle\ =\ 0.
\]

\begin{corollary}  The pairing of gradients $\langle\gradl_{\beta},\gradl_{\beta'}\rangle$ is continuous in a neighborhood of a point $p$ of $\caT(\sigma)\subset\Tbar$ for $\beta,\beta'$ disjoint from the simplex $\sigma$.  The matrix of pairings for a combined short and relative length basis $\basis$ determines a germ at $p$ of a continuous map of $\Tbar$ into a real linear group $GL(\mathbb R)$.
\end{corollary}

Explicit refined counterparts of Masur's expansion already appeared in \cite[Sec. 3]{DW2}, \cite[Theorem 2]{Wlcomp} and \cite[Corollary 4.1]{LSY1}.  The referenced expansions follow a different approach and begin with Masur's holomorphic expansion in $(s,t)$ coordinates for a family of holomorphic $2$-differentials and involve estimating the family hyperbolic metric.  Detailed estimates for the modulus of continuity of the pairings $\langle\gradl_{\beta},\gradl_{\beta'}\rangle$ are provided by Otal \cite[Lemmas 26, 29]{Otwp}.  Otal combines Riera's formula, estimates for the double coset sums, constructions of quasiconformal maps and hyperbolic geometry.  

The expansions reveal the {\em asymptotic factorization} of the WP metric. Near a point $p$ of the augmentation $\caT(\sigma)$ the metric is almost the product of metrics for the normal complex lines $\{\lambda_{\alpha},J\lambda_{\alpha}\}$ and the tangential metric for $T\caT(\sigma)$.   Lower-dimensional WP metrics combine to provide the metric for $T\caT(\sigma)$.  The asymptotic factorization will also appear in the connection and curvature tensor.   The Hessian of geodesic-length is also interesting \cite{Wlbhv}.
\begin{theorem}
\label{hess}
The complex and real Hessians of geodesic-length are uniformly comparable
\[
\partial\overline{\partial}\ell_{\alpha}\ \le\ \hess\ell_{\alpha}\ \le\ \partial\overline{\partial}\ell_{\alpha}.
\]
The first and second derivatives of geodesic-length satisfy
\[
0\ <\ 2\ell_{\alpha}\hess\ell_{\alpha}\,-\,(d\ell_{\alpha})^2\,-\,3(d\ell_{\alpha}\circ J)^2 \quad\mbox{is}\quad O(\ell_{\alpha}^3\|\ \|_{WP}^2)
\]
and
\[
0\ <\ \ell_{\alpha}\partial\overline{\partial}\ell_{\alpha}\,-\,2\partial\ell_{\alpha}\overline{\partial}\ell_{\alpha}\quad\mbox{is}\quad O(\ell_{\alpha}^3\|\ \|_{WP}^2)
\]
where for $c_0$ positive the remainder term constants are uniform for $\ell_{\alpha}\le c_0$.
\end{theorem}

Noting that the Riemannian Hessian is directly related to covariant differentiation $\hess h(U,V)=\langle D_U\grad h,V\rangle$ for vector fields $U,V$ we have the following.

\begin{theorem}
\label{conn}
The covariant derivative of the root-length gradient satisfies 
\[
D_U\lambda_{\alpha}=3\ell_{\alpha}^{-1/2}\langle J\lambda_{\alpha},U\rangle J\lambda_{\alpha}\ +\ O(\ell_{\alpha}^{3/2}\|U\|_{WP})
\]
where for $c_0$ positive the remainder term constant is uniform for $\ell_{\alpha}\le c_0$.  
\end{theorem}

The asymptotic factorization for the connection is as follows \cite{Wlext}.    
\begin{theorem}
The covariant derivatives for the length gradients for $\alpha,\alpha'\in\sigma$, $\beta,\beta'\in\tau$ satisfy: \\
$D_{J\lambda_{\alpha}}\lambda_{\alpha}=3/(2\pi\ell_{\alpha}^{1/2})J\lambda_{\alpha}+O(\ell_{\alpha}^{3/2}),$\;
$D_{\lambda_{\alpha'}}\lambda_{\alpha}=O(\ell_{\alpha}^{3/2})$,\;\\   
$D_{J\lambda_{\alpha'}}\lambda_{\alpha}=O(\ell_{\alpha}(\ell_{\alpha'}^{3/2}+\ell_{\alpha}^{1/2}))$, for 
$\alpha\ne\alpha'$,\;  
$D_{\grad\ell_{\beta}}\lambda_{\alpha}=O(\ell_{\alpha})$, \\ 
$D_{\lambda_{\alpha}}\grad\ell_{\beta}=O(\ell_{\alpha}^{1/2})$,\; $D_{J\lambda_{\alpha}}\grad\ell_{\beta}= 
O(\ell_{\alpha}^{1/2})$,\; and \; $D_{\grad\ell_{\beta'}}\grad\ell_{\beta}$ 
is continuous at $\caT(\sigma)$  with value zero if $\beta$ and $\beta'$ lie on distinct parts.
The remainder term constants depend only on the constant for short geodesic-length.
\end{theorem}

The only covariant derivatives with non trivial leading terms are for differentiation in the normal complex lines and tangential to $\caT(\sigma)$.  Remaining evaluations of derivatives follow from the observation that the almost complex structure $J$ is parallel.  Liu, Sun and Yau also provide a general formula for the covariant derivative \cite{LSY1}.  In \cite{Wlext} we use a quantitative form of the asymptotic factorization for $D$ to show that geodesics in $\caT(\sigma)$ are $C^1$-approximated by geodesics in $\caT$ with tangent field having small projection on $\{\lambda_{\alpha},J \lambda_{\alpha}\}$.  As application we show that the distance between horocycles on surfaces is a convex function on $\caT$. 

A formula for the curvature tensor was developed in \cite{Schu, Siuwp, Wlchern}.  Estimates have since been developed in a number of works \cite{Zh,Zh2,Zh3,Trap} with an at length investigation initiated in \cite{LSY1, LSY2}.  We apply our techniques from \cite{Wlbhv} to understand the curvature tensor evaluated on the root-length gradients.   In the works \cite{Rier, Wlbhv} the Petersson product $\int\varphi\overline{\psi}(ds^2)^{-1}$ was taken as giving the Riemannian cometric; in accordance with Bochner \cite[formula (3)]{Boch} the Hermitian form includes a factor of $\frac12$.  In accordance with Bochner we consider the curvature tensor $R(v,u,z,w)$ as complex linear in $u,z$ and complex anti-linear in $v,w$.  The behavior of the curvature tensor for root-length gradients $\lambda_*=\grad\ell_*$ is straightforward.  
\begin{theorem}
\label{curv}
The curvature tensor evaluated on root-length gradients satisfies
\[
R(\lambda_{\alpha},\lambda_{\alpha},\lambda_{\alpha},\lambda_{\alpha})= \frac{3}{16\pi^3\ell_{\alpha}}\ +\ O(\ell_{\alpha})
\]
and for $\alpha,\beta,\gamma,\delta$ disjoint, not all the same
\[
R(\lambda_{\alpha},\lambda_{\beta},\lambda_{\gamma},\lambda_{\delta}) = O((\ell_{\alpha}\ell_{\beta}\ell_{\gamma}\ell_{\delta})^{1/2})
\]   
where for $c_0$ positive the remainder term constants are uniform for lengths bounded by $c_0$.
\end{theorem}

Liu, Sun and Yau develop estimates for the curvature tensor and an expansion for the Ricci curvature \cite[Lemmas 4.3, 4.5; displays below (4.8); Corollary 4.2]{LSY1} and most importantly apply their approach to analyze the second derivatives of the Ricci tensor. The authors' first applications (see below) begin with a study of the negative of the Ricci form considered as a K\"{a}hler metric.  The resulting metric is shown to have negative holomorphic sectional curvature and an expansion is provided for the resulting Ricci form.  The authors show that the resulting metric is uniformly comparable to important domain metrics, including the asymptotic Poincar\'{e} metric and the K\"{a}hler-Einstein metric \cite{LSY1,LSY2}.   A comparison of the present and the authors' expansions is straightforward with the relations $\ell_{\alpha}\approx 2\pi^2(-\log |t_{\alpha}|)^{-1}$ and $2^{1/2}\pi^2\lambda_{\alpha}\approx (-\log |t_{\alpha}|)^{3/2}t_{\alpha}\frac{\partial}{\partial t_{\alpha}}$.  

The above result determines the curvature tensor since for a pants decomposition the gradients $\{\lambda_{\alpha}\}_{\alpha\in\caP}$ provide a global frame for $T\caT$ over $\mathbb C$.   The asymptotic factorization  applies to the curvature tensor as follows: the leading term of $R(\lambda_{\alpha},\lambda_{\alpha},\lambda_{\alpha},\lambda_{\alpha})$ depends only on $\ell_{\alpha}$; $R(\lambda_{\alpha},\lambda_{\alpha'},\lambda_{\gamma},\lambda_{\delta})$ is at most $O((\ell_{\alpha}\ell_{\alpha'})^{1/2})$; $R(\lambda_{\alpha},\lambda_{\beta},\lambda_{\gamma},\lambda_{\delta})$ is at most $O(\ell_{\alpha}^{1/2})$ and $R(\lambda_{\beta},\lambda_{\beta'},\lambda_{\beta''},\lambda_{\beta'''})$ is continuous at $\caT(\sigma)$. Near a stratum the curvature tensor is almost block diagonalized according to the product of the normal complex lines and the tangent space of the stratum.  The block entries are predicted by asymptotic factorization.

McMullen introduced a K\"{a}hler hyperbolic metric with K\"{a}hler form $\approx\omega_{WP}\,+\,ic\sum_{\alpha\in\caP}\partial\overline{\partial}\operatorname{Log} \ell_{\alpha}$ and established a uniform comparison to the Teichm\"{u}ller metric \cite{McM}. K\"{a}hler hyperbolicity provides a form of average negative holomorphic sectional curvature with exponential volume growth for the properly embedded complex submanifolds of $\caT$.  Applications include the orbifold Euler characteristic having sign alternate with the parity of the genus \cite{Grkh, McM}.  In \cite{LSY1} Liu, Sun and Yau considered the comparison between the WP K\"{a}hler and Ricci forms in their Theorem 6.1 and a uniform comparison between the negative Ricci form and McMullen's metric in their Theorem 6.3.  We combine Theorem \ref{wppair} on the pairing with Theorem \ref{hess} on the Hessian with the above expansion for the curvature tensor to obtain the following refinements of the Liu, Sun and Yau comparisons
\[
\partial\overline{\partial}\log \ell_{\alpha}\,=\,\frac{\partial\ell_{\alpha}}{\ell_{\alpha}}\frac{\overline{\partial}\ell_{\alpha}}{\ell_{\alpha}}\,+\,O(\ell_{\alpha}\|\ \|_{WP}^2)\,=\,\frac{|\langle\lambda_{\alpha},\ \rangle|^2}{\ell_{\alpha}}\,+\, O(\ell_{\alpha}\|\ \|_{WP}^2)
\]
and
\[
R(\lambda_{\alpha},\lambda_{\alpha},\ ,\ )\,=\, \,\frac{3|\langle\lambda_{\alpha},\ \rangle|^2}{4\pi\ell_{\alpha}}\,+\, O(\|\ \|_{WP}^2).
\]

The expansions for the WP metric, connection and curvature tensor hold uniformly on Bers regions for a suitable choice of constant $c_0$.  With Theorem \ref{wppair} precise local bounds for sectional curvatures follow.  Sectional curvatures are uniformly bounded for pairs of basis elements $(\lambda_{\alpha},(J)\lambda_{\alpha'}),\ \alpha,\alpha'\in\sigma$ except for the choices $(\lambda_{\alpha},J\lambda_{\alpha})$.  Sectional curvatures are overall bounded below by $\min_{\alpha\in\sigma}\frac{-3-\epsilon}{\pi\ell_{\alpha}}$ and can be as large as $\min_{\alpha,\alpha'\in\sigma, \alpha\ne\alpha'}\ell_{\alpha}\ell_{\alpha'}$.  The first bound is more specific than and the second bound is larger than Huang's corresponding estimate \cite{Zh,Zh2}.  Huang and Teo consider the genus dependence of bounds for sectional, Ricci and scalar curvatures on the subset of surfaces with injectivity radii bounded below \cite{Zh3,LPT}. 

\section{Methods}
The principal quantities of interest are the gradient $\grad \ell_{\alpha}$ and the operator $(D-2)^{-1}$ for $D$ the hyperbolic Laplacian.  We consider harmonic Beltrami differentials as the Kodaira-Spencer representatives of infinitesimal deformations and represent quantities on the upper half plane $\mathbb H$.  The upper half plane provides the advantages of a single model space and a coordinate model for hyperbolic geometry.  The quantities under investigation are given as integrals of sums over Uniformization group elements of translates of elementary functions.

The WP curvature tensor is given in terms of the integral
\[
\int_R \mu_{\alpha}\overline{\mu_{\beta}}\, (D-2)^{-1} \mu_{\gamma}\overline{\mu_{\delta}} \,dA
\]
for harmonic Beltrami differentials $\mu=\overline{\varphi}(ds^2)^{-1}$, $\varphi\in Q(R)$, and $dA$ the hyperbolic area element, see \cite{Wlchern} and also \cite{Schu, Siuwp}.   The integral defines a $4$-tensor for the tangent space of $\caT$.  To better understand the integral we consider the explicit harmonic Beltrami differentials $\mu_{\alpha}=\grad\ell_{\alpha}$ where $\mu_{\alpha}=\overline{\Theta_{\alpha}}(ds^2)^{-1}$ for $\Theta_{\alpha}\in Q(R)$.   In particular for a geodesic $\alpha$ which lifts to the imaginary axis in $\mathbb H$, the Uniformization group $\Gamma$, and $\Gamma(\alpha)\subset\Gamma$ the stabilizer of the imaginary axis then the Petersson series is
\[
\Theta_{\alpha}\,=\,\sum_{A\in \Gamma(\alpha)\backslash\Gamma}\,A^*\bigl((\frac{dz}{z})^2\bigr).
\]

An essential consideration is the analysis of quantities on collars about geodesics.  Our approach is to introduce Fourier expansions, treat contributions of zeroth terms explicitly and then estimate remaining contributions.  A component in $\mathbb H$ of the lift of the collar about $\alpha$ is the sector $\{\ell_{\alpha}\le \theta \le \pi-\ell_{\alpha}\}$ for the variable $z=re^{i\theta}$ on the upper half plane.  We find for $\ell_{\alpha}$ small the expansions
\begin{displaymath}
\mu_{\alpha} = \left\{ \begin{array}{ll}
a_{\alpha}\sin ^2\theta\,e^{2i\theta} + O((e^{-\theta/\ell_{\alpha}}+e^{-(\pi-\theta)/\ell_{\alpha}})\sin^2\theta )  \quad \textrm{on the collar}\\
\\
\quad O(\ell_{\alpha}^2) \quad\quad\quad\quad\textrm{on the collar complement}
\end{array} \right.
\end{displaymath}
and $a_{\alpha}=\frac{2}{\pi}\,+\,O(\ell_{\alpha}^3)$.  The description of $\mu_{\alpha}$ enables explicit bounds.  The special form of the first remainder provides refined bounds in integrals.  The remainder estimate results from an annular Schwarz lemma applied to $\Theta_{\alpha}$ as follows.  A holomorphic function $h$ on a concentric annulus $\{r_0\le |w|\le r_1\}$ is naturally the sum of a constant $h_0$, a function $h_+$ holomorphic in $\{|w|\le r_1\}$ vanishing at the origin, and a function $h_-$ holomorphic in $\{r_0\le |w|\}$ vanishing at infinity.  From the Cauchy integral formula away from the boundaries the functions $h_-$, $h_0$ and $h_+$ are bounded in terms of an overall bound for $h$.  An overall bound is provided by considering the maximum of $h$ for neighborhoods of the boundaries.  The Schwarz lemma then provides bounds $|h_+|\le |w|/r_1\max_{|w|=r_1-\epsilon}|h_+|$ and $|h_-|\le r_0/|w|\max_{|w|=r_0+\epsilon}|h_-|$.  A change of variable provides the bound for an annulus or collar uniformized by $\mathbb H$.  For holomorphic quadratic differentials evaluation of the constant $h_0=a_{\alpha}$ is given by the pairing with $\Theta_{\alpha}$ \cite{Gardtheta}.  For geodesic-length gradients the constants are given by Theorem \ref{wppair}.  

The following expansion is an essential ingredient towards understanding the main integral
\[
-(D-2)^{-1}|\mu_{\alpha}|^2(z) = \left\{ \begin{array}{ll}
\frac{|a_{\alpha}|^2}{4}\sin ^2\theta\ +\  O(\ell_{\alpha}^2) & \textrm{on the collar where } z=re^{i\theta}\\
\\
\quad O(\ell_{\alpha}\,e^{-dist(\alpha,z)}) & \textrm{on the collar complement.}
\end{array} \right.
\]
The principal term is obtained from the first expansion and the elementary formula $(D-2)\sin^2\theta\,=\,-4\sin^4\theta$.

The techniques for estimating sums over Uniformization group elements are illustrated by considering a basic estimate for the integral kernel $G_2(z,z_0)$ for the operator $(D-2)^{-1}$.  The kernel is given as a sum $G_2(z,z_0)\,=\,\sum_{A\in\Gamma}Q_2(z,Az_0)$.  The associated Legendre function $Q_2$ is negative with a logarithmic pole at zero and for large distance $\delta(z,z_0)$ on $\mathbb H$ satisfies $-Q_2(z,z_0)\le c e^{-2\delta(z,z_0)}$.  The corresponding estimate for the geodesic-length gradient is $|\frac{dz}{z}|^2\sin^2\theta \le c e^{-2\delta(z,\operatorname{axis}_{\alpha})}$.  A basic estimate for the Greens function for $R=\mathbb H/\Gamma$ is as follows \cite[Appendix A.4]{Wlhyp}.
\begin{lemma}
For $\delta(z,z_0)\ge 1$, and $\operatorname{inj}$ injectivity radius and $\delta_R$ distance on $R$ the Greens function satisfies
\[
-Q_2(z,z_0)\,\le\, -G_2(z,z_0)\,\le\, c\,\min_{\zeta=z,z_0}\operatorname{inj}(\zeta)^{-1} \,e^{-\delta_R(z,z_0)}.
\]
\end{lemma}
\noindent{\bf Proof.}\ The lower bound is immediate.  We observe from the collar or Margulis lemma that the overlap number for the action of $\Gamma$ on $\mathbb H$ is bounded by the inverse of the injectivity radius.  In particular for $B(z_0)$ a unit ball in $\mathbb H$ the number of group translates $A(B(z_0))$, $A\in\Gamma$, intersecting $B(z_0)$ is bounded by $c \operatorname{inj}(z_0)^{-1}$.  Next we observe that $-G_2$ (and more generally eigenfunctions of $D$) satisfies a mean value inequality
\[
-G_2(z,z_0)\,\le\, c\int_{B(z_0)}-G_2(z,\zeta)\,dA.
\]
We continue with the rhs as follows
\begin{multline*}
=\,\int_{B(z_0)}\sum_{A\in\Gamma}-Q_2(z,A\zeta)\,dA\,=\,\sum_{A\in\Gamma}\int_{A(B(z_0))}-Q_2(z,\zeta)\,dA\,\\ \le\,c \operatorname{inj}(z_0)^{-1}\int_{\cup_{A\in\Gamma}A(B(z_0))}-Q_2(z,\zeta)\,dA\,\\ \le\,c \operatorname{inj}(z_0)^{-1}\int_{\{\delta(z,\zeta)\,\ge\,\delta_R(z,z_0)-1\}} e^{-2\delta(z,\zeta)}\,dA 
\end{multline*}         
where the first inequality follows since the inverse injectivity radius bounds the overlap number for the covering; the second inequality follows from the general bound for $-Q_2$ and since $\delta_R(z,z_0)$ is a lower bound for the distance of $z$ to $\Gamma z_0$.  The upper bound now follows from the formula $dA\,=\,\sinh\delta\,d\delta d\theta$ and integration.  The argument is complete. 

A region of small (approximately $<1/2$) injectivity radius on a hyperbolic surface is either a collar (a neighborhood of width $2\log 2/\ell_{\alpha}$ about a short geodesic $\alpha$) or a cusp region (a unit area neighborhood of a cusp).  An understanding of small injectivity radius begins with the observation.
\begin{lemma}
For collars and cusp regions $\mathcal C$ the product $\operatorname{inj}(z)e^{\delta(z,\partial\mathcal C)}$ is bounded above and below by universal positive constants.
\end{lemma}

In place of analyzing the operator $(D-2)^{-1}$ Liu, Sun and Yau make use of Schumacher's equation $(\square +1)\,v_j\cdot v_{\bar k}\,=\,A_j\cdot A_{\bar k}$ for the complex Laplacian, horizontal lifts $v_*$ of tangent vectors, and the harmonic Beltrami differentials $A_*$; see \cite[Lemma 2.8]{Schu} and \cite[Sec. 3, opening discussion and Lemma 3.1]{LSY1}.  In brief, considering equations $(D-2)f=h$ is replaced by considering equations $\overline{\partial}f=h$.

\section{Applications of curvature}

\subsection{The work of Liu, Sun and Yau}
Liu, Sun and Yau are investigating questions in the algebraic geometry of the moduli space through their extensive work on curvature \cite{LSY1,LSY2,LSY3,LSY4,LSY5}.  The authors consider the canonical metrics on $\caT$ and $\mathcal M$ with particular focus on the K\"{a}hler-Einstein metric.  They begin by considering the negative of the WP Ricci form as a K\"{a}hler metric.  Degeneration expansions are developed for the first four derivatives of the WP metric.  In Corollary 4.2 of \cite{LSY1} an expansion for the WP Ricci form is provided and in Theorem 4.4 an expansion is provided for the sectional curvatures of the Ricci metric.  The Ricci metric is found to have holomorphic sectional, bisectional and Ricci curvatures bounded from above and below. In an effort to obtain a metric with pinched negative sectional curvature they introduce a perturbation by adding a multiple of the WP metric to the negative Ricci form.  For appropriate multiples the perturbed metric has holomorphic sectional and Ricci curvatures with negative upper and lower bounds, important properties for their considerations.  The work of Farb and Brock \cite{BF} provides that the perturbed metrics can neither be Gromov hyperbolic nor have general sectional curvatures with negative upper and lower bounds.  

The authors compare the canonical metrics.  They combine their expansions with information about the complex geometry of Teichm\"{u}ller space, the Schwarz lemma of Yau, properties of McMullen's K\"{a}hler hyperbolic metric $\approx\omega_{WP}\,+\,ic\sum_{\alpha}\partial\overline{\partial}\operatorname{Log} \ell_{\alpha}$, and the work of Cheng-Mok-Yau on the K\"{a}hler-Einstein metric.   Their results include the following. 
 
\begin{theorem}  The following metrics for $\caT$ and $\mathcal M$ are uniformly comparable: negative Ricci, perturbed negative Ricci, asymptotic Poincar\'{e}, McMullen, K\"{a}hler-Einstein, Teichm\"{u}ller-Kobayashi, Bergman and Carath\'{e}odory.
\end{theorem}

The canonical metrics with distinctive properties lie in a single equivalence class.  They also study the K\"{a}hler-Ricci flow from the Ricci metric to the K\"{a}hler-Einstein metric and apply the methods of Yau to obtain higher order estimates for the canonical metric. 

\begin{theorem}
The K\"{a}hler-Einstein metric has bounded geometry with covariant derivatives of curvature all uniformly bounded.
\end{theorem}

The metrics are singular in a neighborhood of the compactification divisor $\mathcal D$ in the Deligne-Mumford compactification.  The authors establish that Chern forms represent the expected cohomology classes.

\begin{theorem}  Let $\overline{E}=T^*_{\overline{\mathcal M}}(-\log \mathcal D )$ be the logarithmic cotangent bundle of $\overline{\mathcal M}$ and $E$ its restriction to $\mathcal M$.  The WP, Ricci and perturbed Ricci metrics are Mumford good.  The K\"{a}hler-Einstein metric satisfies a generalized Mumford good condition.  The associated Chern form $c_1(\overline{E})$ is positive and $\overline{E}$ is Mumford stable with respect to $c_1(\overline{E})$.
\end{theorem}

The authors also establish the following.
\begin{theorem}
\label{dualNak}
The WP metric for $T_{\mathcal M}$ is dual Nakano negative.
\end{theorem}

Schumacher considers holomorphic families of canonically polarized compact complex manifolds \cite{Schucurv2}.  He finds that the Hermitian metric for the relative canonical bundles induced by the K\"{a}hler-Einstein metrics on the fibers is semi-positive. As an application he finds that the WP metric itself is Nakano positive.  Applications of Theorem \ref{dualNak} include that the Chern numbers of the logarithmic cotangent bundle are positive and a general equivalence of the $L^2$ cohomology for the Ricci, perturbed Ricci and K\"{a}hler-Einstein metrics with the corresponding sheaf cohomology.  An application is the general infinitesimal rigidity of the pair $(\Mbar,\mathcal D)$.  An application with Ji is a Gauss-Bonnet theorem which combined with the computation of Zagier gives the following. 
\begin{theorem}  For the WP, Ricci, perturbed Ricci and K\"{a}hler-Einstein metrics the top Chern classes satisfy
\[
\int_{\mathcal M_g} c_{3g-3} \ =\ \chi(\mathcal M_g)\ =\ \frac{B_g}{4g(g-1)}
\]
for the orbifold Euler characteristic and appropriate Bernoulli number.   
\end{theorem}

\subsection{The model metric $4dr^2\,+\,r^6d\theta^2$}

Yamada initiated the WP synthetic geometry and the comparison of the WP metric for normal complex lines to a standard model metric
\[
(\langle\lambda_{\alpha},\lambda_{\alpha}\rangle)^{-1}(\langle\ ,\lambda_{\alpha}\rangle^2\,+\,\langle\ ,J\lambda_{\alpha}\rangle^2)\,=\,\pi^3(4dr_{\alpha}^2\,+\,r_{\alpha}^6d\theta_{\alpha}^2)\,+\,O(\ell_{\alpha}^3\langle\ ,\ \rangle)
\]
for $\alpha\in\sigma$.  The metric comparison begins with the assignment $2\pi^2 r_{\alpha}^2\,=\,\ell_{\alpha}$ and $\frac{\partial\,}{\partial r_{\alpha}}$ to $2^{3/2}\pi^2\lambda_{\alpha}$.  Model expansions are provided in the references \cite{DW2,Otwp, Wlcomp, Yam2}.  There are various candidate definitions for the angle parameter $\theta_{\alpha}$ in the literature \cite{Abbook,DW2,Msext,Wlbhv}; definitions include the FN angles for a decomposition of hyperbolic surfaces into pairs of pants, as well as the argument of the plumbing parameter for the plumbing family $zw=t$.  Each definition involves auxiliary choices: for the first a pants decomposition is required and for the second attaching maps are required for inserting the plumbing fixture.  We introduce a canonical notion of angle variation as follows.  The counterpart of $\frac{\partial\,}{\partial\theta}$ is taken as the normalized twist $\mathbf T_{\alpha}=(2\pi)^{-1}\ell_{\alpha}^{3/2}J\lambda_{\alpha}=(2\pi)^{-1}\ell_{\alpha}t_{\alpha}$ (recall the relation $2t_{\alpha}=J\gradl_{\alpha}$ between twist and length gradient).  The trajectories of the vector field $\mathbf T_{\alpha}$ are exactly the orbits of the FN deformation with a $2\pi$ flow exactly a full Dehn twist.  We use metric duality to measure the increment in angle.
\begin{definition}
The FN gauge $1$-form is $\rho_{\alpha}\,=\,2\pi(\ell^{3/2}_{\alpha}\langle\lambda_{\alpha},\lambda_{\alpha}\rangle)^{-1}\langle\ , J\lambda_{\alpha}\rangle$.
\end{definition}
The FN gauge satisfies $\rho_{\alpha}(\mathbf T_{\alpha})=1$, $\rho_{\alpha}(\mathbf T_{\alpha'})=O(\ell_{\alpha'}^{3/2})$ for $\alpha\ne\alpha'$, and the definition does not involve auxiliary choices.   We expect that the gauge is not a closed $1$-form (see \cite[Lemma 4.16]{Wlbhv} for further properties and bounds on $d\rho_{\alpha}$).  The model metric expansion is also valid with FN gauges.  

We describe the geometry of the model metric.  The comparison accounts for the normal line expansion for the WP metric, connection and curvature tensor.  The $2$-dimensional model metric is non complete on $\mathbb R^2-\{0\}$.  Coordinate vector fields are evaluated with $\langle\dr,\dr\rangle=4\pi^3,\ \langle\dr,\dthet\rangle=0$ and $\langle\dthet,\dthet\rangle=\pi^3r^6$.  The Riemannian connection $D$ is determined through the relations for vector fields $D_U\langle V,W\rangle=\langle D_UV,W\rangle+\langle V, D_UW\rangle$ and $D_UV-D_VU=[U,V]$.   The formulas for the model connection are
\[
D_{\dr}\tfrac{\partial}{\partial r}=0,\quad D_{\dthet}\tfrac{\partial}{\partial r}=D_{\dr}\tfrac{\partial}{\partial\vartheta}=\tfrac3r\tfrac{\partial}{\partial\vartheta}\quad\mbox{and}\quad D_{\dthet}\tfrac{\partial}{\partial\vartheta}=\tfrac{-3}{4}r^5\tfrac{\partial}{\partial r}.
\]
The formulas are combined to evaluate the curvature tensor $R(U,V)W=D_UD_VW-D_VD_UW-D_{[U,V]}W$.  The Lie bracket of coordinate vector fields $[\dr,\dthet]$ vanishes and from the above evaluations 
$R(\dr,\dthet)\dthet=\frac{-3}{2}r^4\dr$.  The $2$-plane $\dr\wedge\dthet$ has area $2\pi^3r^3$ and the sectional curvature of the model metric is
\[
\frac{\langle R(\dr,\dthet)\dthet,\dr\rangle}{\|\dr\wedge\dthet\|^2}=\frac{-3}{2\pi^3r^2}. 
\]

Theorem \ref{wppair} provides for the WP metric
\[
2\pi\langle \lambda_{\alpha},\lambda_{\alpha}\rangle\,=\,1\ +\ O(\ell_{\alpha}^3)
\]
and Theorem \ref{conn} for the connection
\[
D_{J\lambda_{\alpha}}\lambda_{\alpha}\,=\, \frac{3}{2\pi\ell_{\alpha}^{1/2}}J\lambda_{\alpha}\, +\, O(\ell_{\alpha}^{3/2})
\]
and following Bochner Theorem \ref{curv} for the sectional curvature 
\[
\frac{R(\lambda_{\alpha},\lambda_{\alpha},\lambda_{\alpha},\lambda_{\alpha})}{\langle\lambda_{\alpha},\lambda_{\alpha}\rangle^2}\,=\,\frac{-3}{\pi\ell_{\alpha}}\,+\,O(\ell_{\alpha})\,=\,\frac{-3}{2\pi^3 r_{\alpha}^2}\,+\, O(r_{\alpha}^2).
\]

We relate for the span of $\{\lambda_{\alpha}\}_{\alpha\in\sigma}$ the model metric $\sum_{\alpha\in\sigma} \pi^3(4dr_{\alpha}^2+r_{\alpha}^6\rho_{\alpha}^2)$ to the WP pairing $\langle\ ,\ \rangle$.  The comparison is based on $\frac{\partial\ }{\partial r_{\alpha}}$ having as analog $2^{3/2}\pi^2\lambda_{\alpha}$ and $\frac{\partial\ }{\partial \vartheta_{\alpha}}$ having as analog the Fenchel-Nielsen angle variation $\mathbf T_{\alpha}=(2\pi)^{-1}\ell_{\alpha}^{3/2}J\lambda_{\alpha}$.  The K\"{a}hler form for the model metric is $\sum_{\alpha\in\sigma}2\pi^3r^3_{\alpha}dr_{\alpha}d\theta_{\alpha}=(4\pi)^{-1}\sum_{\alpha\in\sigma}\ell_{\alpha}d\ell_{\alpha}d\theta_{\alpha}$ which for the Fenchel-Nielsen angle $\vartheta_{\alpha}=2\pi\tau_{\alpha}/\ell_{\alpha}$ corresponds to the established formula $\frac12\sum_{\alpha\in\sigma}d\ell_{\alpha}d\tau_{\alpha}$, 
\cite{Wldtau, Wlcusps}.  Theorem \ref{wppair} provides equality of norms modulo the higher order $O$-term. The above expansion for the pairing 
$\langle\lambda_{\alpha},\lambda_{\alpha}\rangle=\frac{1}{2\pi}+O(\ell_{\alpha}^3)$,  Theorem \ref{conn} and the definitions for $2\pi^2r_{\alpha}^2=\ell_{\alpha}$ and $\mathbf T_{\alpha}$ combine to provide the following covariant derivative formulas
\begin{gather*}
D_{\lambda_{\alpha}}\lambda_{\alpha}=O(\ell_{\alpha}^{3/2}),\quad D_{\mathbf T_{\alpha}}(2^{3/2}\pi^2\lambda_{\alpha})=D_{2^{3/2}\pi^2\lambda_{\alpha}}\mathbf T_{\alpha}=\frac{3}{r_{\alpha}}\mathbf T_{\alpha}+O(\ell_{\alpha}^3)\\ \mbox{and}\quad
D_{\mathbf T_{\alpha}}\mathbf T_{\alpha}=\frac{-3}{4}r_{\alpha}^5(2^{3/2}\pi^2\lambda_{\alpha})+O(\ell_{\alpha}^{9/2})
\end{gather*}
in direct correspondence to the formulas for the model connection. The model metric and WP metric for normal lines $C^2$-compare with a factor of $1\,+\,O(r_{\alpha}^4)$. 

The model metric also has a representation as a surface of revolution $(g(u),h(u)\cos v, h(u)\sin v)$ for the pair of functions $g\,=\,\pi^{3/2}\int (4\,-\,9r ^4)^{1/2}dr$ and $h\,=\, \pi^{3/2}r^3$.  The surface of revolution description and also the almost orthonormal tangent frame $\{(2\pi)^{1/2}\lambda_{\alpha},(2\pi)^{1/2}J\lambda_{\alpha}\}$ provide a $u$-Clairaut parameterization with $\lambda_{\alpha}$ defining the meridian direction \cite{doC}.  Two general properties of WP geodesics correspond to asymptotic factorization and the surface of revolution description.  First the component in the complex line $\{\lambda_{\alpha},J\lambda_{\alpha}\}$ of the tangent field of a geodesic should almost be constant length.  In fact from Theorem \ref{conn} for a geodesic $\zeta(t)$ it follows that
\[
\dt (\langle\zeta'(t),\lambda_{\alpha}\rangle^2\,+\,\langle\zeta'(t),J\lambda_{\alpha}\rangle^2)\,=\,O(\ell_{\alpha}^{3/2}(\zeta(t))).
\]
The product of pairings $\langle\zeta'(t),\lambda_{\alpha}\rangle/(\langle\zeta'(t),\zeta'(t)\rangle\langle\lambda_{\alpha},\lambda_{\alpha}\rangle)^{1/2}$ gives the cosine of the angle between $\zeta'(t)$ and the meridian.  With the above given relation the corresponding product of pairings of $\zeta'(t)$ and $J\lambda_{\alpha}$ has the formal interpretation as the sine of the angle. From asymptotic factorization and the classical Clairaut relation the product of the radius $r_{\alpha}^3$ and the sine of the meridian angle should be almost constant. In fact from Theorem \ref{conn} it follows that 
\[
\dt\,\ell_{\alpha}^{3/2}\langle\zeta'(t),J\lambda_{\alpha}\rangle\,=\,O(\ell_{\alpha}^3(\zeta(t))).
\]
  
\subsection{Projection and distance to a stratum}

A complete convex subset of a $CAT(0)$ space is the base for an orthogonal projection \cite[Chap. II.2]{BH}.  Example complete convex subsets of $\Tbar$ include axes of appropriate elements of $MCG$ and sublevel sets of geodesic-length functions (since glf's are convex).  Immediate examples of the latter are the closures of individual strata of $\Tbar$.   Correspondingly for a simplex $\sigma$ of $C(F)$ there are the closure of the stratum $\overline{\caT(\sigma)}$, the orthogonal projection 
$\Pi_{\sigma}$ to $\overline{\caT(\sigma)}$ and the distance $d_{\sigma}$ to the stratum.  General properties include that $d_{\sigma}$ is a convex function and the fibers of $\Pi_{\sigma}$ are fibered by geodesics, called {\em projecting geodesics}.

There are direct relationships between $d_{\sigma}$ and the geodesic-length functions $\ell_{\alpha},\, \alpha\in\sigma$, as well as between the tangent field of a projecting geodesic and the root-length gradients $\lambda_{\alpha},\,\alpha\in\sigma$ in \cite{Wlbhv}.  We found expansions for distance
\[
d_{\sigma}\,\le\,(2\pi\sum_{\alpha\in\sigma}\ell_{\alpha})^{1/2} \quad\mbox{and}\quad d_{\sigma}\,=\,(2\pi\sum_{\alpha\in\sigma}\ell_{\alpha})^{1/2}\ + \ O(\sum_{\alpha\in\sigma}\ell_{\alpha}^{5/2}).
\]
Theorem \ref{conn}  provides for a refined expansion for the tangent field of a unit speed projecting geodesic
\[
\dt \,=\, (2\pi)^{-1}\sum_{\alpha\in\sigma}a_{\alpha}\lambda_{\alpha}\ +\ O(t^4)
\]
for small $t$ with $(2\pi)^{1/2}\|(a_{\alpha})\|_{Euclid}=1$.   The expansions for distance and the tangent field manifest the asymptotic factorization.  Projecting geodesics are to $C^1$ fourth-order approximated by constant sums of $\lambda_{\alpha},\,\alpha\in\sigma$.  Projecting geodesics are almost geodesics of the metric for the product of $\{\lambda_{\alpha},J\lambda_{\alpha}\}$ and $T\caT(\sigma)$. 

The first and second variations of the projection distance $d_{\sigma}$ are interesting for applications.  We explain now that for top dimensional strata the second variation of distance $d_{\sigma}$ is given by the classical second variation formula for Jacobi fields.

Differentiability is considered as follows.  For $\sigma=\{\alpha\}$, we write $\caT(\alpha)$, $\Pi_{\alpha}$ and $d_{\alpha}$ for the associated quantities.  Given $\epsilon$ positive, $d_{\alpha}$ is approximated by the projection distance $d_{\alpha,\epsilon}$ to the complete convex sublevel set $\{\ell_{\alpha}\le\epsilon\}$.  The first and second variation of  $d_{\alpha,\epsilon}$ are given by the classical formulas.   Provided the first and second variations have the corresponding expressions for $d_{\alpha}$ as their uniform limits, then differentiability and application of the formulas for $d_{\alpha}$ follow. The first variation is indeed given by the classical cosine formula \cite[Chapter II.3, Corollary 3.6]{BH}.  An important step is showing that Jacobi fields for projecting geodesics to $\{\ell_{\alpha}=\epsilon\}$ converge to Jacobi fields of projecting geodesics to $\caT(\alpha)$.  The fibers of the projection $\Pi_{\alpha}$ are fibered by projecting geodesics.  We now describe the Jacobi fields that describe the variation through projecting geodesics of the projecting geodesic $\eta(t)$ with $\eta(0)\in\caT(\alpha)$, $\eta$ with outer endpoint $p\in\caT$ and tangent field $\dt$.  The local frame $\{\lambda_{\alpha},J\lambda_{\alpha},\gradl_{\beta}\}_{\beta\in\tau}$ provides a frame at $\eta(0)$.  
\begin{lemma}
\label{jac}
Given the projecting unit-speed geodesic $\eta(t)$ to $\caT(\alpha)$ and a vector $v$ at the endpoint $p$ in $\caT$ there exists a unique Jacobi field $V$ with $V_p=v$ and initial $C^2$-expansion 
\[
V(t)\,=\,at\dt\ +\ bt^3J\dt\ +\ A\ +\ O(t^4)
\]
where $a\|\eta\|=\langle v, \dt\rangle_p$, $A$ is parallel along $\eta$ with $A\perp\dt,\,A\perp J\dt$ and $\|A(0)\|\,<\,\|v\,-\,\langle v,\dt\rangle\dt\|_p$.  The remainder constant is locally bounded in terms of $\eta,\, a,\, b$ and $A$.
\end{lemma}
Consistent with asymptotic factorization the first component of the Jacobi field describes displacement of the endpoint along $\eta$, the second component describes rotation in the infinitesimal normal complex line ($\frac{d}{d\theta}=r^3J\frac{d}{dr}$ is infinitesimal rotation for the model metric) and the third component describes displacement orthogonal to $\dt,\,J\dt$ in $\caT(\alpha)$. 

The length $\|\eta\|$ of the projecting geodesic is the distance $d_{\alpha}(p)$.  For $p$ close to the stratum $\caT(\alpha)$ we note that $a\|\eta\|\,=\,\langle v,\dt\rangle$, $b\|\eta\|^3\,=\,\langle v,J\dt\rangle\,+\, O(\|\eta\|^4|v|)$ and the parallel field $A$ satisfies $A\perp\dt,\,A\perp J\dt$ with the vector expansion at the endpoint
\[
v\,=\,a\|\eta\|\dt\,+\,b\|\eta\|^3J\dt\,+\,A\,+\, O(\|\eta\|^4|v|).
\]

The approximation of the distance $d_{\alpha}$ by $d_{\alpha,\epsilon}$ provides a setting for deriving the second variation formula.  The classical formula \cite[Chap. 1, \S 6]{CE} is combined with Theorem \ref{conn}, Lemma \ref{jac} and the approximation to show for Jacobi fields $V,W$ as above that the second variation of distance is
\[
\ddot{d_{\alpha}}[V,W]\, =\, \langle D_WV,\dt\rangle_p\,+\,\langle D_{\dt}(V-a_V\dt),(W-a_W\dt)\rangle_p.
\]
The derivative $D_WV$ vanishes for the variation of distance along a geodesic.


\end{document}